\input ./style/arxiv-general.cfg
\documentclass[MSNbibl,number,seceqn,citesort,dvips]{arxbj}
\makeatletter
   \@ifpackageloaded{graphicx}{}{\usepackage{graphicx}}
\makeatother
\usepackage{mathbh}

\volume{22}
\issue{2}
\pubyear{2016}
\firstpage{653}
\lastpage{680}
\doi{10.3150/14-BEJ641}
\docsubty{FLA}

\makeatletter
\newcommand{\ds}{\displaystyle}
\newtheorem{theorem}{Theorem}

\newtheorem{proposition}{Proposition}
\makeatother

\begin{document}
\begin{frontmatter}

\title{The number of accessible paths in the~hypercube}
\runtitle{The number of accessible paths in the hypercube}

\begin{aug}
\author[A]{\inits{J.}\fnms{Julien}~\snm{Berestycki}\corref{}\thanksref{e1}\ead[label=e1,mark]{julien.berestycki@upmc.fr}},
\author[B]{\inits{\'{E}.}\fnms{\'{E}ric}~\snm{Brunet}\ead[label=e2]{Eric.Brunet@lps.ens.fr}}
\and
\author[A]{\inits{Z.}\fnms{Zhan}~\snm{Shi}\ead[label=e3,mark]{zhan.shi@upmc.fr}}
\address[A]{Laboratoire de Probabilit\'es et Mod\`eles Al\'eatoires, UPMC,
75005 Paris, France.\\ \printead{e1,e3}}
\address[B]{Laboratoire de Physique Statistique, ENS, UPMC, CNRS, 75005 Paris,
France.\\ \printead{e2}}
\end{aug}

%
\received{\smonth{11} \syear{2013}}

%
\begin{abstract}
Motivated by an evolutionary biology question,
we study the following problem: we consider the hypercube $\{0,1\}^L$ where
each node carries an independent random variable uniformly distributed on
$[0,1]$, except $(1,1,\ldots,1)$ which carries the value $1$ and
$(0,0,\ldots,0)$ which carries the value $x \in[0,1]$. We study the number
$\Theta$ of paths from vertex $(0,0,\ldots,0)$ to the
opposite vertex $(1,1,\ldots,1)$ along which the values on the nodes
form an increasing sequence. We
show that if the value on $(0,0,\ldots,0)$ is set to $x=X/L$ then
$\Theta/L$ converges in
law as $L \to\infty$ to $\mathrm{e}^{-X}$ times the product of two standard
independent exponential variables.

As a first step in the analysis, we study the same question when the graph
is that of a tree where the root has arity $L$, each node at level 1 has
arity $L-1$, \ldots, and the
nodes at level $L-1$ have only one offspring which are the leaves of the
tree (all the leaves are assigned the value~1, the root the value~$x\in[0,1]$).
\end{abstract}

%
\end{frontmatter}

\section{Introduction and motivation}

\subsection{The model}

We consider the following problem: for $L\ge 1$, let $(x_\sigma,
\sigma\in\{0,1\}^L)$ be a sequence of i.i.d. random variables with
uniform distribution on $[0,1]$ except for $\sigma_0=(0,0,\ldots,0)$
for which we fix $x_{\sigma_0}=x$ with $x$ given and $\sigma
_L=(1,1,\ldots,1)$ for which we fix $x_{\sigma_L}=1$.
Viewing $\{0,1\}^L$ as the $L$-dimensional hypercube we ask how many
oriented paths there are from $\sigma_0$ to
$\sigma_L$:
\[
\sigma_0 \to\sigma_1 \to\sigma_2\to\cdots
\to\sigma_L,
\]
where each $\sigma_{i+1}$ is obtained from $\sigma_i$ by changing a single
0 into a 1 in the sequence $\sigma_i$, such that values
$x_\sigma$ form an increasing sequence:
\[
x_{\sigma_0}< x_{\sigma_1}<\cdots< x_{\sigma_L}.
\]
Such paths are said to be \textit{open} or \textit{accessible}.
A variant of this model which we also consider is when the value
$x_{\sigma_0}$ at the starting point is picked randomly as the other
$x_\sigma$.

\subsection{Motivation}

This question is motivated by some classical and recent works in
evolutionary biology.
Consider the following very simplified model for the evolution of an
organism. The genetic information of the organism is encoded into its genome
which, for our purposes is a chain of $L_0$ sites. With time,
the organism accumulates mutations which are only single site substitutions.

If we suppose that there are only two possible alleles on each site, it
makes sense when looking at the genome to only record
whether the allele carried at a given site is the
original one (the ``wild type'') or the mutant. We will represent a
genetic type by a sequence
of 0's and 1's of length $L_0$ where we put a 0 at position $i$ if this site
carries the original code or a 1 if it carries the mutant. Hence, a genetic
type is a point $\sigma\in\{0,1\}^{L_0}$, the $L_0$-dimensional hypercube.
This is a classical way of encoding the possible evolutionary states of a
population. For instance, in their seminal article \cite{kauffmann-levin}
Kauffman and Levin write:

\begin{quote}
Consider as a concrete example, a space of $N$ peptides constrained
to use two amino acids, say leucine and alanine. Then 1 and 0 can
represent the two amino acids, and each peptide is a binary string,
length $N$ of 1 and 0 values. Such strings are easily represented as
vertices of an $N$-dimensional Boolean hypercube.
\end{quote}

As an organism evolves by successive mutations, its genetic type travels
along the edges of the hypercube. Each genetic type $\sigma
\in\{0,1\}^{L_0}$ is characterized by a certain fitness value
$x_{\sigma}$.

Assume that the population is in a regime with a low
mutation rate and strong selection; this means that when a new genetic
type (mutant) appears in a resident population, it must either fixate
(i.e.,  it
invades the whole population and becomes the resident type) if it has
better fitness
or become extinct (i.e., no one in the population carries
this type after some time) if its fitness is lower. Therefore, in that
low mutation and strong
selection regime, the only possible evolutionary paths are such that the
fitness is always increasing. We say that such paths are \emph{open}.
In biology, paths with increasing fitness values are also referred to as
selectively accessible (see \cite{weinreich-watson-chao,weinreich-delaney-depristo-hartl,franke-kloezer-devisser-krug}).
The idea that the population moves as a whole along the vertices of the
hypercube is also classical and can be found in \cite{kauffmann-levin}
or in \cite{klg}:

\begin{quote}
One can think of the adaptive process as a continuous time,
discrete state Markov process, in which the entire
population is resident at one state and then jumps with fixed
probabilities to each of its 1-step mutant fitter variants.
\end{quote}

Somewhere in the $L_0$-dimensional hypercube, there is a type with the
highest fitness. We call~$L$ the distance between that type and the
original one; that is,  $L$ is the number of mutated alleles in the fittest
type. A natural question is whether there is an open path from the
original type
to the fittest. Such a path has at least $L$ steps but may contain many
more. Because of the low mutation rate, evolution takes time, and we
interest ourselves here only in the \emph{shortest} open paths leading to
the fittest type, that is in the paths with exactly $L$ steps, for which
mutation never goes back: a~site can only change from the original type to
the mutant one.

In that setting, it is thus sufficient to consider the $L$-dimensional hypercube
which contains (as opposing nodes) the original type, noted
$\sigma_0=(0,0,\ldots,0)$ and the fittest type $\sigma_L=(1,1,\ldots
,1)$. We
consider paths through that hypercube along the edges which always move
further away from the origin (i.e.,  at each step a 0 is changed into
a 1 in the sequence) and, out of the $L!$ possible paths, we wish to count
the number~$\Theta$ of open paths, that is the number of paths such that
the fitness values form an increasing sequence which represent a possible
evolutionary history of the organism towards its optimum fitness.

To count the number of open paths, we need to have a model for the fitness
values $x_\sigma$ of
all the nodes. There are many choices as to how to do this. One could, for
instance, choose $x_\sigma$ to be the number of ones in $\sigma$. In that
case, all direct paths from $\sigma_0$ to $\sigma_L$ are accessible. However
this corresponds to a very smooth, linear {\it fitness landscape} which
does not match observations.

Instead, our choice is to pick the fitness values as independent random
variables
with a common
distribution. As we are only interested in whether a sequence is
increasing or not, the results will not depend on the specific distribution
(as long as it is atomless). We therefore choose to give a
fitness~$x_{\sigma_L}=1$ to the
fittest node and to assign uniform random numbers between 0 and 1 to each
other node. This is the so-called ``House of Cards'' model, which was
introduced by Kingman \cite{kingman}, which is also \cite{altenberg} the
$\mathit{NK}$ model studied by \cite{kauffmann-levin} in the limit $K=N-1$.
The variant we consider where $x_{\sigma_0}$ is not randomly picked
but fixed
to a given value $x$ has also been considered recently \cite{carneirohartl,klozner}.
In a follow-up paper \cite{BBZ2}, we explore the situation in which we
authorize arbitrary paths on the hypercube, and not only the shortest ones.

\subsection{A toy model}
\begin{figure}[t]

\includegraphics{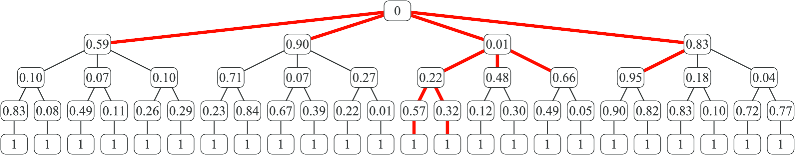}

\caption{A tree for $L=4$ and $x=0$.
The bold lines (also red in the online version) indicate
the
directed paths going down from the root (at the top) and which visit an
increasing sequence of numbers. There are $\Theta=2$ paths going all
the way
down to the leaves of the tree.}
\label{figuretree}
\end{figure}

The correlation structure of the hypercube raises significant technical
challenges. As a first step, we study the following simplified problem:
instead of working on the $L$-dimensional hypercube we chose to work on
a deterministic rooted tree as in Figure~\ref{figuretree} with arity
decreasing from $L$ to~1: the root is connected to $L$ first level nodes,
each first level node is connected to $L-1$ second level nodes, etc. There
are $L$ levels in the tree and $L!$ directed paths. The number of possible
steps at level $k$ is then $L-k$, as on the hypercube. Each of the $L!$
leaves of the tree (at level $L$) are assigned the value~1. All the other
nodes are assigned independent random numbers uniformly drawn between~0
and~1, except perhaps the root to which we may choose to give a fixed
value~$x$. We are interested in the number $\Theta$ of directed paths on
the tree going from the root to one of the leaves where the numbers
assigned to the visited nodes form an increasing sequence. As before, such
a path is said to be open.

%

We mention that other models of paths with increasing fitness values on
trees have also been considered in the literature \cite{nowak-krug,RZ2013,xinxin}.

\section{Main results}

We call
\[
\Theta= \mbox{number of open paths on the tree or the hypercube}.
\]
We want to describe the law of $\Theta$.
Whether we work on the tree or on the hypercube will always be made
clear from
the context. The starting position for a path on the tree is the root
of the tree; by abuse
of language we call ``root of the hypercube'' the starting position
$\sigma_0=(0,0,\ldots,0)$ of a path on the hypercube. Throughout the paper, we use
the following notations for the probability of
an event and for the expectation and the variance of a number:
\begin{eqnarray}
&& \mathbb{P}^x(\cdot),\qquad \mathbb{E}^x(\cdot),\qquad
\operatorname{Var}^x(\cdot) \qquad\mbox{when the root has value $x$},
\nonumber\\
&& \mathbb{P}^*(\cdot)=\int_0^1\mathrm{d} x
\mathbb{P}^x(\cdot),\qquad \mathbb{E}^*(\cdot)=\int
_0^1\mathrm{d} x \mathbb{E}^x(
\cdot), \qquad\operatorname{Var}^*(\cdot)\nonumber\\
\eqntext{\mbox{when the root is uniform
in~$[0,1]$.}}
\end{eqnarray}
Note that the size $L$ of the tree or of the hypercube is implicit in
the notation.
The notation $a_L \sim b_L$,
$L\to\infty$ (which we write in words as ``$a_L$ is equivalent to
$b_L$ when $L$ is large'') means $\lim_{L\to\infty} \frac{a_L}{b_L} =1$.

Obtaining the expectation of the number of open paths when the starting value
has the fixed value $x$ is easy:
there are $L!$ paths in the tree or the hypercube. Each path has probability
$(1-x)^{L-1}$ that the $L-1$ intermediate numbers between the root and
the leaf are between $x$ and $1$. Furthermore, there is probability
$1/(L-1)!$ that these intermediate numbers form an increasing sequence.
Hence, the probability that a given path is open is
$(1-x)^{L-1}/(L-1)!$ and thus
%
\begin{equation}
\label{1stmoment}
\mathbb{E}^x(\Theta) = L (1-x)^{L-1}
\end{equation}
both for the tree and the hypercube.

Thus, if $x=0$ there are on average $L$ open paths, and if $x>0$ the
number of open paths goes in probability to zero when $L$ becomes large.

The most biologically relevant variant of the model is when $x$ is also
randomly picked as the other nodes. The expectation of $\Theta$ (both
on the tree and on the hypercube) is trivially obtained by integrating
(\ref{1stmoment}):
%
\begin{equation}
\mathbb{E}^*(\Theta)=1,
\end{equation}
but the typical number of paths for $L$ large is not of order 1, as can
be seen by looking at the variance of $\Theta$:
%
\begin{equation}
\label{Evstar}
\lim_{L\to\infty}\frac{\operatorname{Var}^*(\Theta)}L=1\qquad \mbox{on the
tree}.
\end{equation}
(All the variance computations on the tree are carried out in
Section~\ref{sec2nd}.) In fact, this can be understood by considering
starting values~$x$ scaling with the size $L$ of the system as $x=X/L$ with
$X\ge0$ fixed:

\begin{proposition}\label{Pprop1}
In the case of the tree,
%
\begin{equation}
\label{EvL}
\lim_{L\to\infty}\mathbb{E}^{X/L}(\Theta/L) =
\mathrm{e}^{-X}, \qquad \lim_{L\to\infty}\operatorname{Var}^{{X}/L}(
\Theta/L) = \mathrm{e}^{-2X}.
\end{equation}
In the case of the hypercube (Hegarty--Martinsson \cite{hegarty-martinsson})
%
\begin{equation}\label{varhyper}
\lim_{L\to\infty}\mathbb{E}^{{X}/L}(\Theta/L) =
\mathrm{e}^{-X}, \qquad \lim_{L\to\infty}\operatorname{Var}^{{X}/L}(
\Theta/L) = 3\mathrm{e}^{-2X}.
\end{equation}
\end{proposition}

(Note that Hegarty--Martinsson consider a different scaling regime, but
their proof can be adapted without any modification to the result above.)

From Proposition~\ref{Pprop1} for $x=X/L$ and $L$ large, the variance
of $\Theta$ scales
like the square of the expectation of $\Theta$. This means that when the
starting value $x$ is $\mathcal{O}(1/L)$, the number $\Theta$ of open
paths is $\mathcal{O}(L)$, like its expectation. When $x$ is
chosen randomly, there is a probability~$\mathcal{O}(1/L)$ that
$x=\mathcal{O}(1/L)$ yielding $\mathcal{O}(L)$ open
paths. On average we thus expect
$\mathcal{O}(1)$ open paths, with a variance $\mathcal{O}(L)$, as in
(\ref{Evstar}).

This heuristic can be made more precise.

\begin{theorem}\label{thmtree}
On the tree, for a starting value $x=X/L$ (with $X\ge0$ fixed), the
variable $\Theta/L$ converges in
law when $L\to\infty$ to $\mathrm{e}^{-X}$ multiplied by a standard
exponential variable.
\end{theorem}

\begin{theorem}\label{conj}
On the hypercube, for a starting value $x=X/L$ (with $X\ge0$ fixed),
the variable $\Theta/L$ converges in
law when $L\to\infty$ to $\mathrm{e}^{-X}$ multiplied by the product of two
independent standard exponential variables.
\end{theorem}

It will become apparent in the proofs that we get a product of two
independent variables on the hypercube because, locally near both corners
$(0,0,\ldots, 0)$ and $(1,1,\ldots,1)$, the hypercube graph looks roughly
like the tree.

We conclude with the following remark: when the starting value $x$ is
picked randomly, even if the expectation and the variance of $\Theta$ is
dominated by values of $x = \mathcal{O}(1/L)$, the probability that
there exists at
least one open path is dominated by starting values $x = (\ln L)/L+\mathcal{O}(1/L)$.
This was made clear on the hypercube in \cite{hegarty-martinsson} and we
here state the tree counterpart.

\begin{theorem}\label{thmtreep0}
On the tree, when the starting value is $x= (\ln L+X)/L$
%
\begin{equation}
\label{EnlogL}
\lim_{L\to\infty} \mathbb{E}^{({\ln L+X})/L} (\Theta)=
\mathrm{e}^{-X}, \qquad\lim_{L\to\infty} \operatorname{Var}^{({\ln L+X})/L}
(\Theta)= \mathrm{e}^{-2X}+\mathrm{e}^{-X}.
\end{equation}
When the starting value $x$ is chosen at random uniformly in
$[0,1]$, the probability to have no open path goes to 1 as $L\to\infty$ and
%
\begin{equation}\label{limit}
\mathbb{P}^*(\Theta\ge1)\sim\frac{\ln L}L \qquad \mbox{as }L\to\infty.
\end{equation}
\end{theorem}

The rest of the paper is organized as follows: we start by proving
Theorem~\ref{thmtreep0} in Section~\ref{sproofthm2} as it is the
simplest, and then we prove Theorem~\ref{thmtree} in
Section~\ref{secproofthm1}. The proofs rely on Proposition~\ref{Pprop1}
which is itself proven in Section~\ref{sec2nd} for the tree.
In Section~\ref{seccascade}, we introduce the notion of Poisson cascade,
which allows to give a probabilistic interpretation of one of the main
objects introduced in our proofs. Finally, Theorem~\ref{conj} (on the
hypercube) is proven in Section~\ref{secconj}.

\section{Proof of Theorem \texorpdfstring{\protect\ref{thmtreep0}}{3}}
\label{sproofthm2}

In (\ref{EnlogL}), the result on the expectation is trivial from
(\ref{1stmoment}), and the result on the variance is obtained in
Section~\ref{sec2nd}. In this section, we prove~(\ref{limit}).

Let us start with the upper bound. Markov's inequality with
(\ref{1stmoment}) leads to
%
\begin{equation}
\mathbb{P}^x(\Theta\ge1)\le\min \bigl[1, L(1-x)^{L-1}
\bigr].
\end{equation}
We split the integral $\mathbb{P}^*(\Theta\ge1) = \int_0^1 \mathbb
{P}^x(\Theta\ge
1)\, \mathrm{d}x$ at $x_0 =1-\exp[-(\ln L)/(L-1)]$ since that is the point
such that $L(1-x_0)^{L-1}=1$. We end up with
%
\begin{equation}
\label{E6}
\mathbb{P}^*(\Theta\ge1)\le 1 - \exp \biggl[-\frac{\ln L }{L-1}
\biggr]+\exp \biggl[-\frac{L}{L-1}\ln L \biggr] =\frac{\ln L}L+\mathcal{O}
\biggl(\frac{1}L \biggr).
\end{equation}

We now turn to the lower bound.
Let $L\mapsto f(L)$ be a function diverging more slowly than $\ln L$:
%
\begin{equation}
\lim_{L\to\infty} f(L)=\infty,\qquad  0\le f(L) \le\ln L, \qquad \lim
_{L\to\infty}\frac{f(L)}{\ln L}=0.
\end{equation}
It is sufficient to show that
%
\begin{equation}\label{condcor}
\lim_{L\to\infty}\mathbb{P}^{({\ln L-f(L)})/L}(\Theta\ge1)=1,
\end{equation}
because
%
\begin{equation}\label{lb}
\mathbb{P}^*(\Theta\ge1)\ge\int_0^{({\ln L-f(L)})/L}\mathbb
{P}^x(\Theta\ge1) \,\mathrm{d} x\ge\frac{\ln L-f(L)}L
\mathbb{P}^{({\ln L-f(L)})/L}(\Theta\ge1),
\end{equation}
where we used that $x\mapsto\mathbb{P}^x(\Theta\ge1)$ is a
non-increasing function.
Taking $L$ large in (\ref{lb}) assuming (\ref{condcor}) gives the
lower bound.

It now remains to show (\ref{condcor}).
We consider a tree started from $x=[\ln L-f(L)]/L$, and call~$m$ the number
of nodes at the first level with a value between $x$ and $(\ln L)/L$.
One has $m \sim\operatorname{Bin}(L, f(L)/L)$, that is $m$ is a binomial of
parameters $L$ and $f(L)/L$. Conditionally on $m$,
the probability to have no open path in the tree is smaller than the
probability to have no open path through these $m$ specific nodes. Thus by
summing over all possible values of $m$, we get:
%
\begin{eqnarray}
&& \mathbb{P}^{({\ln L-f(L)})/L}(\Theta=0)
\nonumber
\\[-8pt]
\label{binom}
\\[-8pt]
\nonumber
 && \quad\le\sum_{m=0}^L
\pmatrix{L \cr m} \biggl(\frac{f(L)}L \biggr)^m \biggl(1-
\frac{f(L)}L \biggr)^{L-m} \bigl[\mathbb{P}^{({\ln L})/L;L-1}(\Theta=0)
\bigr]^m,
\end{eqnarray}
where we used $\mathbb{P}^x(\Theta=0)\le\mathbb{P}^{(\ln
L)/L}(\Theta=0)$ for
$x\le(\ln L)/L$. Note the obvious extension to the notation to mark that
the probability on the right-hand side is for a tree of size~$L-1$ and not
$L$ as on the left-hand side. Summing~(\ref{binom}), one gets
%
\begin{equation}\label{binom2}
\mathbb{P}^{({\ln L-f(L)})/L}(\Theta=0)\le \biggl[1-\frac{f(L)}L \bigl(1-
\mathbb{P}^{({\ln L})/L;L-1}(\Theta =0) \bigr) \biggr]^L.
\end{equation}
But from Cauchy--Schwarz (applied to $\Theta$ and $\mathbh{1}(\Theta
\ge1)$) and
(\ref{EnlogL}), which is proved in Section~\ref{proofn2}, one has
%
\begin{equation}
\mathbb{P}^{({\ln L})/L}(\Theta\ge1)\ge\frac{\mathbb{E}^{({\ln L})/L}(\Theta)^2}{\mathbb{E}^{({\ln L})/L}(\Theta^2)}\mathop{\longrightarrow}_{L\to \infty}\frac{1}3,
\end{equation}
so that, for $L$ large enough
%
\begin{equation}
\mathbb{P}^{({\ln L})/L;L-1}(\Theta\ge1) \ge\mathbb{P}^{({\ln(L-1)})/({L-1});L-1}(\Theta\ge1)
\ge0.33,
\end{equation}
and thus, for $L$ large enough
%
\begin{equation}
\mathbb{P}^{({\ln L-f(L)})/L}(\Theta= 0)\le \biggl[1-\frac{f(L)}L
\mathbb{P}^{({\ln L})/L;L-1}(\Theta\ge 1) \biggr]^L \le \biggl[1 - 0.33
\frac{f(L)}L \biggr]^L
\end{equation}
which goes to zero as $L\to\infty$, as required.

\section{Proof of Theorem \texorpdfstring{\protect\ref{thmtree}}{1}}
\label{secproofthm1}

In this section, we consider the case of the tree with a starting value $x$
which scales as $x=X/L$, $X\ge0$ being a fixed number. The natural starting
point for a proof would be to introduce $G$, the generating function
of~$\Theta$:
%
\begin{equation}\label{gen0}
G(\lambda,x,L)=\mathbb{E}^x \bigl(\mathrm{e}^{-\lambda\Theta} \bigr),
\end{equation}
with parameter $\lambda\ge0$, for which it is very easy to show from the
tree geometry (each of the $L$ nodes at the first level is the root of an
independent tree of size~$L-1$) that:
%
\begin{equation}\label{rec0}
G(\lambda,x,1)=\mathrm{e}^{-\lambda},\qquad  G(\lambda,x,L)= \biggl[x+\int
_x^1 \mathrm{d}y G(\lambda,y,L-1)
\biggr]^L \qquad \mbox{for }L>1.
\end{equation}
However, extracting the limiting distribution directly from (\ref{rec0})
seems difficult because the number of levels and the size of
each level increase together and because the fixed point equation does not
give the $\lambda$ dependence of the result. We shall rather use an idea
which proved to be very generic and powerful in branching processes: the
value of a random variable is decided during the early stages of
a branching process; at later stages the law of large
numbers kicks in (see, e.g., \cite{lalley-sellke}).

Assume that all the information at the first $k$ levels of the tree is
known, and call~$\Theta_k$ the expected number of paths given that information:
%
\begin{equation}\label{thetadefinition}
\Theta_k=\mathbb{E}(\Theta|\mathcal{F}_k),
\end{equation}
where $\mathcal{F}_k$ is the available information up to level~$k$.
For instance, consider the tree of Figure~\ref{figuretree} up to
level~$k=2$. There are three paths still open with end values (at
level~2) given by $0.22$, $0.66$ and $0.95$. Therefore from
(\ref{1stmoment}), $\Theta
_2=2(1-0.22)+2(1-0.48)+2(1-0.66)+2(1-0.95)=3.38$. Similarly,
$\Theta_1 = 3(1-0.59)^2+3(1-0.90)^2+3(1-0.01)^2+3(1-0.83)^2 \approx
3.56$ and
$\Theta_3=\Theta_4=\Theta=2$. A general expression of $\Theta_k$
for $k<L$ is
%
\begin{equation}\label{expnk}
\Theta_k=\sum_{|\sigma|=k}
\mathbh{1}_{\{\sigma \ \mathrm{open}\}} (L-k) (1-x_\sigma)^{L-k-1},
\end{equation}
where we sum over all nodes $\sigma$ at level $|\sigma|=k$ in the tree,
$x_\sigma$ is the value of the node $\sigma$ and the event $\{\sigma
\mbox{ open} \}$ is the $\mathcal{F}_{|\sigma|}$-measurable event
that the
path from the root to node $\sigma$ is open.

Heuristically, one expects that $\Theta$ and $\Theta_k$ are of order
$L$. When $k$ is small, $\Theta$ has no
reason to be close to $\Theta_k$. However when $k$ is large, there are many
paths open up to level~$k$ which all contribute to the value of~$\Theta$.
The law of large numbers leads to a small variance, given $\mathcal{F}_k$, of
$\Theta/L$ and $\Theta_k/L$ becomes a good approximation of $\Theta/L$.
The advantage of this approach is that one can take the $L\to\infty$ limit
for a fixed~$k$ (keeping the depth of the tree constant) and then
take the $k\to\infty$ limit.

Our proof consists then in two steps:
\begin{itemize}
\item first we show that
%
\begin{equation}\label{nnnk}
\lim_{L\to\infty}\mathbb{P}^{{X}/L} \biggl(\frac{\Theta} L
\le z \biggr) =\lim_{k\to\infty}\lim_{L\to\infty}
\mathbb{P}^{{X}/L} \biggl(\frac{\Theta_k}L \le z \biggr),
\end{equation}
which means that $\Theta/L$ for a starting point $X/L$ and $L$ large
has the
same distribution as $\Theta_k/L$ with the same starting point for large~$L$
and then for large~$k$;
\item then we make use of a generating function similar to~(\ref{gen0})
to show that the distribution of $\Theta_k/L$ after taking the limits
is given
by an exponential law.
\end{itemize}

\subsection{Proof of (\texorpdfstring{\protect\ref{nnnk}}{4.5})}

Pick $\delta>0$. Observe that
%
\begin{eqnarray}
\mathbb{P} \biggl(\frac{\Theta} L\le z\Big|\mathcal{F}_k \biggr)
&\le& \mathbh{1} \biggl(\frac{\Theta_k}L\le z +\delta \biggr) + \mathbb{P}
\biggl(\frac{|\Theta-\Theta_k|}L\ge\delta\Big|\mathcal {F}_k \biggr),
\nonumber
\\[-8pt]
\label{mainineq}
\\[-8pt]
\nonumber
\mathbb{P} \biggl(\frac{\Theta} L\le z\Big|\mathcal{F}_k
\biggr) &\ge& \mathbh{1} \biggl(\frac{\Theta_k}L\le z - \delta \biggr) -
\mathbb{P} \biggl(\frac{|\Theta-\Theta_k|}L\ge\delta\Big|\mathcal {F}_k\biggr).
\end{eqnarray}
The first inequality follows from the simple remark that it is obviously
true when $\Theta_k/L\le z+\delta$ and that, when $\Theta
_k/L>z+\delta$, it
is necessary to have $(\Theta_k-\Theta)/L\ge\delta$ to get $\Theta
/L\le z$.
Similarly, the lower bound in the next line is trivial when
$\Theta_k/L>z-\delta$ and it is sufficient to have
$(\Theta-\Theta_k)/L<\delta$ when $\Theta_k/L\le z-\delta$.

As $\Theta_k$ is the expectation of $\Theta$ given $\mathcal{F}_k$,
one has from
Chebyshev's inequality:
%
\begin{equation}\label{Cheby}
\mathbb{P} \biggl(\frac{|\Theta-\Theta_k|}L\ge\delta\Big|\mathcal {F}_k
\biggr)\le \frac{\operatorname{Var}(\Theta|\mathcal{F}_k)}{L^2\delta^2}.
\end{equation}
We substitute (\ref{Cheby}) into (\ref{mainineq}) and then take the
expectation
(over $\mathcal{F}_k$):
\begin{eqnarray*}
\mathbb{P}^{{X}/L} \biggl(\frac{\Theta} L\le z
\biggr) &\geq &
\mathbb{P}^{{X}/L} \biggl(\frac{\Theta_k}L\le z -\delta \biggr)-
\frac{\mathbb{E}^{{X}/L} [\operatorname{Var}(\Theta|\mathcal{F}_k) ]}{L^2\delta^2},   \\[-3pt]
\mathbb{P}^{{X}/L} \biggl(\frac{\Theta} L\le z
\biggr)  &\le & \mathbb{P}^{{X}/L} \biggl(\frac{\Theta_k}L\le z +\delta
\biggr)+\frac{\mathbb{E}^{{X}/L} [\operatorname
{Var}(\Theta|\mathcal{F}_k) ]}{L^2\delta^2}.
\end{eqnarray*}
Therefore to show (\ref{nnnk}), it is sufficient to have
%
\begin{equation}\label{limvar}
\lim_{k\to\infty}\limsup_{L\to\infty}\frac{1}{L^2}
\mathbb{E}^{{X}/L} \bigl[\operatorname{Var}(\Theta|\mathcal
{F}_k) \bigr]=0,
\end{equation}
as well as the existence and continuity of the right-hand side limit of
(\ref{nnnk}).

In Section~\ref{proofmajo}, we will show by direct analysis of the
second moment
that
%
\begin{equation}\label{majo}
\lim_{L\to\infty}\frac{1}{L^2}\mathbb{E}^{{X}/L} \bigl[
\operatorname{Var}(\Theta|\mathcal{F}_k) \bigr] =\frac{\mathrm{e}^{-2X}}{2^k},
\end{equation}
which yields (\ref{limvar}).
We  now compute the distribution of $\Theta_k/L$
in the double limit $L\to\infty$ and $k\to\infty$ and, as the
result is a continuous function of $z$, this\vspace*{-6pt}  completes the proof.

\subsection{Distribution of \texorpdfstring{$\Theta_k$}{$Theta_k$}}

Similarly to (\ref{gen0}),
we define $G_k(\lambda,x,L)$ the generating function of $\Theta_k$
for a tree
of size $L$ and a value $x$ at the\vspace*{-2pt} root:
%
\begin{equation}\label{Gkdef}
G_k(\lambda,x,L)=\mathbb{E}^x \bigl(
\mathrm{e}^{-\lambda\Theta_k} \bigr).
\end{equation}
As $\Theta_0=\mathbb{E}^x(\Theta)=L(1-x)^{L-1}$ one\vspace*{-2pt} has
%
\begin{equation}\label{initG}
G_0(\lambda,x,L) = \exp \bigl[{-\lambda L(1-x)^{L-1}}
\bigr],
\end{equation}
and the recursion\vspace*{-4pt} relation
%
\begin{eqnarray}
G_k(\lambda,x,L)&=& \biggl[x+\int_x^1
\mathrm{d}y G_{k-1}(\lambda,y,L-1) \biggr]^L
\nonumber
\\[-9pt]
\label{recG}\\[-9pt]
\nonumber
&=&  \biggl[1-\int
_x^1 \mathrm{d}y \bigl(1-G_{k-1}(
\lambda,y,L-1) \bigr) \biggr]^L,
\end{eqnarray}
to be compared to (\ref{rec0}).

This relation is obtained by decomposing on what happens at the first
splitting. For a node~$\sigma$ connected to the root let
$\Theta_k{(\sigma)}$ be the conditional expectation given $\mathcal{F}_k$
of the number of open paths going through $\sigma$. The
$\{\Theta_k{(\sigma)}\}_{ |\sigma|=1}$ is a collection of $L$ independent
$\mathcal{F}_k$-measurable independent variables,\vspace*{-2pt} hence:
%
\begin{equation}
G_k(\lambda,x,L) = \bigl[ \mathbb{E}^x \bigl(
\mathrm{e}^{-\lambda
\Theta_k{(\sigma)}} \bigr) \bigr]^L,
\end{equation}
where $\sigma$ is a given node in the first generation.

Let us evaluate $\mathbb{E}^x ( \mathrm{e}^{-\lambda\Theta_k{(\sigma
)}} )$. If
$x_\sigma<x$, then $\Theta_k{(\sigma)}=0$ and since this event has
probability $x$ it contribute $x \mathrm{e}^{-\lambda0}=x$ to the expectation.
With a probability $\mathrm{d}y$ for $y\in[x,1]$ the value at the node
lies in $(y, y+ \mathrm{d}y)$ and some paths might go through that
node. The subtree rooted at
$\sigma$ is like the initial tree but of dimension $L-1$ and we want to
evaluate the average number of paths in that subtree given the information
after $k-1$ steps, hence the term in the integral of (\ref{recG}).

The strategy is to take the $L\to\infty$ limit at fixed $k$ in
(\ref{initG}) and (\ref{recG}) after a proper rescaling, then to let
$k\to\infty$. We only consider $\lambda\ge0$; it is sufficient to
characterize the distribution, and it simplifies the arguments below.

\begin{longlist}[Step 1.]
\item[Step 1.] We first show that the following limit exists (for $\mu\ge0$):
%
\begin{equation}\label{limG1}
\forall a,b, \qquad  G_k \biggl(\frac{\mu}{L+a},\frac{X} {L+b}, L
\biggr) \mathop{\longrightarrow}_{L\to\infty}\tilde{G}_k(\mu,X),
\end{equation}
and that the limit satisfies
%
\begin{equation}\label{limG2}
\tilde{G}_k(\mu,X)=\exp \biggl[-\int_X^\infty
\bigl[1-\tilde G_{k-1}(\mu,Y) \bigr]\,\mathrm{d}Y \biggr], \qquad
\tilde{G}_0(\mu,X)=\exp \bigl[-\mu \mathrm{e}^{-X} \bigr].
\end{equation}
From (\ref{initG}), it is obvious that (\ref{limG1}) holds for $k=0$ with
the limit given in (\ref{limG2}). Choosing $k>0$, we assume that
(\ref{limG1}) holds for $G_{k-1}$. Then, after a change of variables
in (\ref{recG}),
\[
G_k \biggl(\frac{\mu}{L+a},\frac{X}{L+b},L \biggr)=
\biggl[1-\frac{1}{L+b}\int_X^{L+b}
\mathrm{d}Y \biggl(1-G_{k-1} \biggl(\frac{\mu}{L+a},
\frac{Y}{L+b},L-1 \biggr) \biggr) \biggr]^L. 
\]
The $G_{k-1}$ on the right-hand side has an $L\to\infty$ limit. From its
definition (\ref{Gkdef}), one has
%
\begin{equation}
1\ge G_k(\lambda,x,L)\ge1-\lambda\mathbb{E}^x(
\Theta_k)=1-\lambda \mathbb{E}^x(\Theta)=1-\lambda
L(1-x)^{L-1}.
\end{equation}
Then, assuming $\mu\ge0$, for all $a$ and $b$, one has for $L$ large
enough (depending on $a$ and~$b$):
%
\begin{equation}
1\ge G_{k-1} \biggl(\frac{\mu}{L+a},\frac{Y}{L+b},L-1 \biggr)
\ge1-2\mu \mathrm{e}^{-Y/2}.
\end{equation}
Thus, from the dominated convergence theorem, we have
%
\begin{equation}
\int_X^{L+b} \mathrm{d}Y \biggl(1-G_{k-1}
\biggl(\frac{\mu}{L+a},\frac{Y}{L+b},L-1 \biggr) \biggr)
\mathop{\longrightarrow}_{L\to\infty} \int_X^\infty
\mathrm{d}Y \bigl(1-\tilde{G}_{k-1}(\mu, Y) \bigr),
\end{equation}
and thus
(\ref{limG1}) holds for $G_k$ with the relation (\ref{limG2}).

\item[Step 2.] The fact that (\ref{limG1}) holds means that
when starting with $x=X/L$,
the random variable $\Theta_k/L$ has a well defined limit as $L$ goes to
infinity, and that the generating function of that limit is $\tilde{G}_k$.
We now use the recurrence (\ref{limG2}) to take the
$k\to\infty$ limit which will show that $\lim_{L\to\infty} \Theta
_k/L$ converges
(when $k\to\infty$) to an exponential variable. This task is greatly
simplified by noticing (by a simple recurrence) that one can write
$\tilde G_k$ as a function of one variable only:
%
\begin{equation}
\tilde{G}_k(\mu,X)=F_k\bigl(\mu \mathrm{e}^{-X}
\bigr)
\end{equation}
with
%
\begin{equation}\label{recF}
F_k(z)=\exp \biggl[-\int_0^z
\frac{1-F_{k-1}(z')}{z'}\,\mathrm{d}z' \biggr],\qquad F_0(z)=
\mathrm{e}^{-z}.
\end{equation}

We shall show that the solution to (\ref{recF}) satisfies
%
\begin{equation}\label{limF}
F_k(z)\mathop{\longrightarrow}_{k\to\infty} \frac{1}{1+z}\qquad
\mbox{for } z>-1,
\end{equation}
which implies that $\lim_{L\to\infty} \Theta_k/L$ converges weakly
when $k\to\infty$
to an exponential distribution of expectation $\mathrm{e}^{-X}$.
Note that we only need to consider $z\ge0$ and, in fact, we proved
(\ref{recF}) only
for $z\ge0$, but (\ref{limF}) holds for the solution to (\ref{recF}) for
$z\in(-1,\infty)$.

Defining $\delta_k(z)$ for $z>-1$ and $z\ne0$ by
%
\begin{equation}
F_k(z)=\frac{1}{1+z}-\frac{z^2}{(1+z)^3} \frac{\delta_k(z)}{2^k},
\end{equation}
it is easy to see that there exists a constant $M$ such that
for all $k$ and all $z >- 1$
%
\begin{equation}\label{recdelta}
0\le\delta_k(z)\le M.
\end{equation}
Indeed, for $k=0$,
%
\begin{equation}
\delta_0(z)=\frac{(1+z)^3}{z^2} \biggl(\frac{1}{1+z}-
\mathrm{e}^{-z} \biggr),
\end{equation}
$\delta_0(z)\ge0$ for $z>-1$ because $\mathrm{e}^{z}\ge1+z$ by convexity.
Furthermore, $\delta_0(z)$ can be defined by continuity at $z=0$,
has a limit in $z=+\infty$ and in $z=-1$ and reaches therefore
a maximum~$M$ on $(-1,\infty)$, which initializes (\ref{recdelta}).

Assuming now (\ref{recdelta}) at order $k-1$, one has
%
\begin{equation}
F_k(z)=\frac{1}{1+z}\exp \biggl[-\int_0^z
\frac{z'}{(1+z')^3} \frac{\delta_{k-1}(z')}{2^{k-1}}\,\mathrm{d}z' \biggr],
\end{equation}
leading to
\[
\frac{1}{1+z}\ge F_k(z) 
\ge
\frac{1}{1+z} \biggl[1-\frac{M}{2^{k-1}}\int_0^z
\frac
{z'}{(1+z')^3} \,\mathrm{d} z' \biggr]=\frac{1}{1+z}
\biggl[1-\frac{M}{2^{k-1}} \frac
{z^2}{2(1+z)^2} \biggr],
\]
which gives (\ref{recdelta}) at order~$k$. Hence, the limit (\ref{limF})
holds. This completes the proof of Theorem~\ref{thmtree}.
\end{longlist}

\section{Results on the second moment for the tree}
\label{sec2nd}

The goal of this section is to prove the second moment results
(\ref{Evstar}), (\ref{EvL}), (\ref{EnlogL}) and (\ref{majo}) which were
used in the proofs of Theorems \ref{thmtree} and \ref{thmtreep0}.

\subsection{Exact expression of the second moment}

The expectation of $\Theta^2$ is the sum, over all pairs of paths, of the
probability that both paths are open. There are $L!^2$ pairs of paths in
the system. For a given pair, the probability that they are both open
depends on the number $q\in\{0,1,2,\ldots, L-2,L\}$ of bonds shared
by the
paths. (Note: two paths cannot have exactly $L-1$ bonds in common.) The
following facts are clear:
\begin{itemize}
\item the number of pairs of paths which coincide all the way ($q=L$)
is $L!$;
\item the probability that ``both'' paths in such a pair are open is
$(1-x)^{L-1}/(L-1)!$;
\item the number of pairs of paths which coincide for $q=0,1,\ldots, L-2$
steps and then branch is $L! (L-q-1)(L-q-1)!$
(Remark: $1\cdot 1!+2 \cdot 2!+3\cdot 3!+\cdots+(L-1)(L-1)!=L!-1$, hence, one recovers that
the total number of pairs of path is $L!^2$.);
\item the probability that both paths in such a pair are open is
%
\begin{equation}
\frac{(1-x)^{2L-q-2}}{(2L-q-2)!}\pmatrix{2L-2q-2
\vspace*{2pt}\cr
L-q-1}.
\end{equation}
Indeed, excluding the starting and end points, there are $2L-q-2$
total different nodes in such a pair of paths. All these nodes must
be larger than $x$, hence the $(1-x)^{2L-q-2}$ term. This is however
not sufficient because the values on the nodes must be correctly ordered.
Out of the $(2L-q-2)!$ possible orderings (see the denominator), the
only good ones are those such that the $q$ smallest terms are in their
correct order
in the shared segment (only one choice), and the
$2L-2q-2$ remaining terms are separated into two well ordered blocks of $L-q-1$
terms, one for each path; the only freedom is to choose which terms go
to which path, hence the binomial coefficient.
\end{itemize}
This leads to
%
\begin{equation}\label{2ndmoment}
\mathbb{E}^x\bigl(\Theta^2\bigr)=\sum
_{q=0}^{L-2}a(L,q) (1-x)^{2L-q-2} +L
(1-x)^{L-1},
\end{equation}
where
%
\begin{equation}
\label{aLq}
a(L,q)=\frac{L!(2L-2q-2)!}{(L-q-2)!(2L-q-2)!}.
\end{equation}
The isolated term in (\ref{2ndmoment}) corresponds to the pairs of
identical paths and is equal to $\mathbb{E}^x(\Theta)$.

\subsection{Estimates and bounds on the $a(L,q)$}

Expanding the factorials in $a(L,q)$, one gets
%
\begin{equation}\label{aqsmall}
a(L,q)=\frac{L^2}{2^q}\frac{ (1-{1}/L) (1-{2}/L)\cdots
(1-({q+1})/L)}{
 (1-({q+2})/({2L}) ) (1-({q+3})/({2L}))\cdots (1-({2q+1})/({2L}))}.
\end{equation}
From this expression, one gets
the following equivalent when $L\to\infty$ and $q\ll\sqrt L$:
%
\begin{equation}\label{aqequiv}
a(L,q)=\frac{L^2}{2^q} \bigl[1+\mathcal{O}\bigl(q^2/L\bigr)
\bigr].
\end{equation}
For $q$ close to $L$, one has the values
$a(L,L-2)=2$, $a(L,L-3)=24/(L+1)$, $a(L,L-4)=360/[(L+1)(L+2)]$, etc.

We want to find a good upper bound for $a(L,q)$. We first show that
$q\mapsto\ln a(L,q)$ is a convex function for $L\ge q+3$. Indeed
%
\begin{equation}
\ln a(L,q)-\ln a(L,q-1)=\ln\frac{(L-q-1)(2L-q-1)}{(2L-2q)(2L-2q-1)}
\end{equation}
so that
%
\begin{eqnarray}
&& \bigl[\ln a(L,q)-\ln a(L,q-1)\bigr]-\bigl[\ln a(L,q-1)-\ln a(L,q-2)\bigr]
\nonumber
\\[-8pt]
\label{lna}
\\[-8pt]
\nonumber
&& \quad=\ln\frac{(L-q-1)(2L-q-1)(2L-2q+2)(2L-2q+1) }{(2L-2q)(2L-2q-1)(L-q)(2L-q)}.
\end{eqnarray}
Since the denominator is clearly positive as soon as $L\ge q+1$
we see that $\ln a(L,q)$ is convex if in (\ref{lna}) the numerator is
bigger than the denominator. This condition leads to
%
\begin{equation}
(L-q)^2(2L-1)-(2L-q-1) (2L-2q+1)\ge0,
\end{equation}
which holds as soon as $L-q\ge3$.

Assume $L\ge12$ so that $a(L,L-3)\le2$ and let
%
\begin{equation}
q_0(L)= \biggl\lceil\frac{\ln( L^2)}{\ln2}+1 \biggr\rceil.
\end{equation}
From (\ref{aqsmall}), one has
%
\begin{equation}
a(L,q)\le\frac{L^2}{2^q} \frac{1}{ (1-({2q+1})/({2L}) )^q}.
\end{equation}
Applying this to $q=q_0(L)$,
one easily gets $a(L,q_0(L))\le1$ for $L$ large enough. (The term
$L^2/2^{q_0}$ is smaller than~$1/2$, and the parenthesis converges to~1.)

By using the convexity of $\ln a(L,q)$, one has
\[
\ln a(L,q)\le \cases{\displaystyle\ln a(L,0)+\frac{q}{q_0(L)} \bigl[\ln a
\bigl(L,q_0(L) \bigr)-\ln a(L,0) \bigr],  & \quad$\mbox{for }0\le q\le
q_0(L)$,\vspace*{3pt}
\cr
\ln2, & \quad $\mbox{for }q_0(L)\le q
\le L-2$.}
\]
But $\ln a (L,q_0(L) )\le0$ so that
%
\begin{equation}
a(L,q)\le\cases{ \ds a(L,0)\exp \biggl[-\frac{\ln a(L,0)}{q_0(L)}q \biggr], & \quad $\mbox{for }0
\le q\le q_0(L)$,\vspace*{3pt}
\cr
2, & \quad$\mbox{for }q_0(L)
\le q\le L-2$.}
\end{equation}
Remark now that $a(L,0)=L(L-1)<L^2$ and $\ln a(L,0)/ q_0(L)\to\ln2$ as
$L\to\infty$. This implies that,
for $L$ large enough, $\ln a(L,0)/ q_0(L)>\ln1.99$ and that
%
\begin{equation}
\label{aqbound}
a(L,q)\le\cases{ L^2 1.99^{-q}, & \quad \mbox{for
}$0\le q\le q_0(L)$,\vspace*{3pt}
\cr
2, & \quad\mbox{for }$q_0(L)\le q\le L-2$.}
\end{equation}

\subsection{Proofs of the limits (\texorpdfstring{\protect\ref{Evstar}}{2.3}),
(\texorpdfstring{\protect\ref{EvL}}{2.4}) and (\texorpdfstring{\protect\ref{EnlogL}}{2.6}) 
of \texorpdfstring{$\operatorname{Var}(\Theta)$}{$\operatorname{Var}(Theta)$}}
\label{proofn2}

The second moment (\ref{2ndmoment}) is written as a sum from~$q=0$ to
$q=L-2$. To
prove the various limits we need, the strategy is always the
same:
\begin{enumerate}[(3)]
\item[(1)]
Split the sum over $q$ into two parts; one going from~0 to $q_0(L)$ and
one going from $q_0(L)+1$ to $L-2$.
\item[(2)] In the first sum, replace
$a(L,q)$ by its equivalent~(\ref{aqequiv}); this is justified with the
dominated convergence theorem, using the bound~(\ref{aqbound}).
\item[(3)]
Show that the second sum does not contribute using the bound (\ref{aqbound}).
\end{enumerate}

\begin{pf*}{Proof of (\protect\ref{Evstar})}
Integrating (\ref{2ndmoment}) over $x$ and using $\mathbb{E}^*(\Theta
)=1$, one gets
\begin{eqnarray*}
\frac{\operatorname{Var}^*(\Theta)}{L} &=& \frac{\mathbb{E}^*(\Theta
^2)-\mathbb{E}^*(\Theta)^2}{L}= \sum_{q=0}^{L-2}
\frac{a(L,q)/L}{2L-q-1} \\
&=& \sum_{q=0}^{q_0(L)}
\frac{a(L,q)/L^2}{2-({q+1})/L} +\sum_{q=q_0(L)+1}^{L-2}
\frac{a(L,q)/L}{2L-q-1}.
\end{eqnarray*}
In the first sum, the running term is equivalent to $2^{-q}/2$ when
$L$ is large and is dominated by $1.99^{-q}$ for~$L$ large enough.
Therefore, this first sum converges to~$1$. In the second sum, the running
term is smaller than $2/L^2$, implying that the whole second sum is
smaller than $2/L$ and thus vanishes in the large~$L$ limit.
\end{pf*}

\begin{pf*}{Proof of (\protect\ref{EvL})}
We divide (\ref{2ndmoment}) by $L^2$, replace $x$ by $X/L$ and split the
sum:
%
\begin{eqnarray}
\frac{\mathbb{E}^{{X}/L}(\Theta^2)}{L^2}&=& \sum_{q=0}^{q_0(L)}
\frac{a(L,q)}{L^2} \biggl(1-{\frac{X} L} \biggr)^{2L-q-2}
\nonumber
\\[-8pt]
\\[-8pt]
\nonumber
&&{}+\sum
_{q=q_0(L)+1}^{L-2}\frac{a(L,q)}{L^2} \biggl(1-{
\frac{X} L} \biggr)^{2L-q-2} +\frac{ (1-{X}/L)^{L-1}}L.
\end{eqnarray}
The running term in the first sum is equivalent to $2^{-q}\mathrm{e}^{-2X}$ and is
dominated by $1.99^{-q}$, therefore the first term converges to $2\mathrm{e}^{-2X}$
as $L\to\infty$. The running term in the second sum is smaller than $2/L^2$
implying that the whole second sum is
smaller than $2/L$ and thus vanishes in the large~$L$ limit.
The isolated term goes also to zero. Therefore, the whole expression
converges to $2\mathrm{e}^{-2X}$ and one recovers the variance in (\ref{EvL}) after
subtracting\vspace*{-2pt} $\mathbb{E}^{X/L}(\Theta/L)^2$.
\end{pf*}

\begin{pf*}{Proof of (\protect\ref{EnlogL})}
We now take $x=(\ln L+X)/L$ and split again the sum in (\ref
{2ndmoment}) into
two parts:
%
\begin{eqnarray}
\mathbb{E}^x\bigl(\Theta^2\bigr)&=& \sum
_{q=0}^{q_0(L)}\frac{a(L,q)}{L^2}\times
L^2(1-x)^{2L-q-2}
\nonumber
\\[-9pt]
\label{split2}
\\[-9pt]
\nonumber
&&{}+\sum_{q=q_0(L)+1}^{L-2}
a(L,q) (1-x)^{2L-q-2}+L(1-x)^{L-1}.
\end{eqnarray}
Using\vspace*{-5pt}
%
\begin{equation}
\lim_{L\to\infty} L^2 \biggl(1-\frac{\ln L+X}L
\biggr)^{2L-q-2}=\mathrm{e}^{-2X},
\end{equation}
into (\ref{split2}), the running term in the first sum is equivalent to
$2^{-q} \mathrm{e}^{-2X}$ and is dominated by $1.99^{-q}(\mathrm{e}^{-2X}+1)$ for $L$ large
enough (because $L^2(1-x)^{2L-q-2}\le L^2(1-x)^{2L-q_0(L)-2}$, which
becomes close to its limit when $L$ gets large). Therefore, the first sum
converges to $2\mathrm{e}^{-2X}$.
We write an upper bound of the second sum of (\ref{split2}) using
$a(L,q)\le2$
and then extending the sum to the interval $[0,L-1]$:
%
\begin{eqnarray}
\sum_{q=q_0(L)+1}^{L-2}a(L,q)
(1-x)^{2L-q-2} &\le &  2(1-x)^{2L-2}\frac{(1-x)^{-L}-1}{(1-x)^{-1}-1}
\nonumber
\\[-9pt]
\\[-9pt]
\nonumber
&\le &  2
\frac{(1-x)^{L-1}}{x} \sim2\frac{\mathrm{e}^{-X}}{\ln L},
\end{eqnarray}
which goes to zero for $L$ large. Finally, the last term in (\ref{split2})
converges to $\mathrm{e}^{-X}$; putting things together, one finds $\mathbb{E}^{(\ln
L+X)/L}(\Theta^2)\to2\mathrm{e}^{-2X}+\mathrm{e}^{-X}$. Removing the expectation
squared, one
recovers\vspace*{-3pt} (\ref{EnlogL}).
\end{pf*}

\subsection{Exact expression for \texorpdfstring{$\mathbb{E}^x[\operatorname{Var}(\Theta|\mathcal{F}_k)]$}{$\mathbb{E}^x[\operatorname{Var}(Theta|\mathcal{F}_k)]$}}

The number $\Theta$ of paths given $\mathcal{F}_k$ is the sum over
all the nodes at
level~$k$ of the number of paths through that node. These variables are
independent; therefore
%
\begin{equation}
\operatorname{Var}(\Theta| \mathcal{F}_k)=\sum
_{|\sigma|=k} \mathbh{1}_{\{\sigma \ \mathrm{open} \}} v(x_\sigma,L-k),
\end{equation}
where $v(x,L)$ is the variance of $\Theta$ for a tree of size
$L$ started at $x$
%
\begin{eqnarray}
v(x,L) &:=& \mathbb{E}^x\bigl(\Theta^2\bigr)-
\mathbb{E}^x(\Theta)^2
\nonumber
\\[-9pt]
\label{5.18}\\[-9pt]
\nonumber
&=&  -L(1-x)^{2L-2}+\sum
_{q=1}^{L-2}a(L,q) (1-x)^{2L-q-2}+L
(1-x)^{L-1}.
\end{eqnarray}

Taking the expectation over $\mathcal{F}_k$, one gets
%
\begin{equation}
\label{5.19}
\mathbb{E}^x\bigl[\operatorname{Var}(\Theta|
\mathcal{F}_k)\bigr]=\frac
{L!}{(L-k)!}\int_x^1
\mathrm{d}x_\sigma \frac{(x_\sigma-x)^{k-1}}{(k-1)!}v(x_\sigma,L-k),
\end{equation}
where $L!/(L-k)!$ is the number of terms in the sum and where the
fraction in the integral is the probability that $\sigma$ is open
given the value
of $x_\sigma>x$.

In (\ref{5.19}), we replace $v$ by its expression given in (\ref
{5.18}), and we integrate each term in the sum to find
\begin{eqnarray*}
E^x\bigl[\operatorname{Var}(\Theta|\mathcal{F}_k)\bigr]
&=& -\frac
{L!(2L-2k)!}{(L-k-1)!(2L-k-2)!}(1-x)^{2L-k-2}
\\
&&{}+\sum_{q=1}^{L-2} a(L-k,q)
\frac
{(2L-2k-q-2)!L!}{(L-k)!(2L-k-q-2)!}(1-x)^{2L-k-q-2}
\\
&&{}+L(1-x)^{L-1},
\end{eqnarray*}
where\vspace*{1.5pt} we have used the formula $\int_x^1 \mathrm{d}z (z-x)^m (1-z)^n
= {m! n! \over(m+n+1)!} (1-x)^{m+n+1}$. Then, after changing $q$ into
$q-k$ and using extensively the expression (\ref{aLq}) of $a(L,q)$ we obtain
\begin{eqnarray}
\mathbb{E}^x\bigl[\operatorname{Var}(\Theta|\mathcal{F}_k)
\bigr]&=& -\frac
{a(L,k)}{L-k-1}(1-x)^{2L-k-2}
\nonumber
\\[-8pt]
\label{Exvarnk}\\[-8pt]
\nonumber
&&{}+\sum_{q=k+1}^{L-2}
a(L,q) (1-x)^{2L-q-2} +L(1-x)^{L-1}.
\end{eqnarray}
Note that apart from the first term, this is exactly the same as the full
variance $v(x,L)$ except that the sum over $q$ begins at $k+1$ instead of
at 1.

\subsection{Proof of (\texorpdfstring{\protect\ref{majo}}{4.9})}
\label{proofmajo}

We now divide (\ref{Exvarnk}) by $L^2$, set $x=X/L$ and consider $L$ large.
We only need an upper bound, but it is as easy to calculate the exact limit.
As in Section~\ref{proofn2}, we split the sum into two parts;
one where the index~$q$ runs from $k+1$ to $q_0(L)$ and one from
$q_0(L)+1$ to
$L-2$. In the first part, using the dominated convergence theorem with the
bound (\ref{aqbound}):
%
\begin{equation}
\lim_{L\to\infty}\sum_{q=k+1}^{q_0(L)}
\frac{a(L,q)}{L^2} \biggl(1-\frac{X} L \biggr)^{2L-q-2}=\sum
_{q=k+1}^\infty \frac{1}{2^q}
\mathrm{e}^{-2X}=\frac{1}{2^k}\mathrm{e}^{-2X}.
\end{equation}
Also using the bound (\ref{aqbound}), the second part of the sum goes
to zero:
%
\begin{equation}
\frac{1}{L^2}\sum_{q=q_0(L)+1}^{L-2}a(L,q)
\biggl(1-\frac{X} L \biggr)^{2L-q-2}\le\frac{1}{L^2} \times L
\times2.
\end{equation}
It is very easy to check that in (\ref{Exvarnk}) the two isolated terms
(divided by $L^2$, of course) go also to zero, so that one finally obtains
(\ref{majo}).

\section{A relation with Poisson cascades}
\label{seccascade}

Our model is closely related to cascades of Poisson processes. In fact, the
arguments we used in Section~\ref{secproofthm1} can be presented in terms
of Poisson cascades. Let us make a brief description.

We recall the sequence of functions $F_k$, $k\ge0$, defined in (\ref{recF}):
%
\begin{equation}
F_k(z) = \exp \biggl[-\int_0^z
\frac{1-F_{k-1}(z')}{z'}\,\mathrm{d}z' \biggr],\qquad  F_0(z)=
\mathrm{e}^{-z}.
\end{equation}
It is clear that $F_0(z)$ is the Laplace transform of the Dirac
measure at 1, and that $F_1$ is the Laplace transform of $\sum_{j=1}^\infty
X_j$, where $(X_j, j\ge1)$ is a Poisson process on $(0, 1]$ with
intensity $\mathbh{1}_{(0, 1]}(x) \frac{\mathrm{d}x} x$.

We now define a cascade of Poisson processes. At generation $k=0$,
there is
only one particle at position $1$. At generation $k=1$, this particle is
replaced by the atoms $(X_j^{(1)}, j\ge1)$ of a Poisson process on $(0,
1]$ with intensity $\mathbh{1}_{(0, 1]}(x)\, \frac{\mathrm{d}x} x$. At
generation $k=2$, for each $j$, the particle at position $X_j^{(1)}$ is
replaced by $(X_j^{(1)}X_{j,\ell}^{(2)}, \ell\ge1)$, where
$(X_{j,\ell}^{(2)}, \ell\ge1)$ is another Poisson process with
intensity $\mathbh{1}_{(0, 1]}(x) \,\frac{\mathrm{d}x}x$ (all the Poisson
processes are assumed to be independent). Iterating the procedure results
in a cascade of Poisson processes. We readily check, by induction on $k$,
that $F_k$ is the Laplace transform of $Y_k$, the sum of the positions at
the $k$th generation of the Poisson cascade.

What was proved in Section~\ref{secproofthm1} can be stated in terms
of the
cascade of Poisson processes.
Recall from (\ref{thetadefinition}) that $\Theta_k=\mathbb{E}(\Theta|\mathcal{F}_k)$.

\begin{theorem}
\label{tcascadepoisson}
\textup{(i)} For any $k\ge0$ and for $x=X/L$, $\frac{\Theta
_k} L$ converges weakly, when
$L\to\infty$, to $\mathrm{e}^{-X}Y_k$.

\textup{(ii)} When $k\to\infty$, $Y_k$ converges weakly to the standard
exponential law.
\end{theorem}

\section{Proof of Theorem~\texorpdfstring{\protect\ref{conj}}{2}}
\label{secconj}

In this section, we adapt the methods used in
Section~\ref{secproofthm1} to obtain the distribution of
the number of open paths on the hypercube when~$L$ goes to infinity.

In the large $L$ limit, both the width (the number of possible moves
at each step) and the depth (the number of steps) on the hypercube go to
infinity, which makes studying the limit difficult.
We worked around that problem on the tree by introducing $\Theta_k$, the
expected number of paths given the information $\mathcal{F}_k$ after
$k$ steps, and
by sending first $L$ (now representing only the width of the tree) and then
$k$ (the depth) to infinity.

We use the same trick on the hypercube, but with a twist:
the hypercube is symmetrical when exchanging the starting and end
points, and there is no reason to privilege one or the other. Therefore,
we call~$\Theta_k$ the expected number of paths in the hypercube given the
information~$\mathcal{F}_k$ at the first $k$ levels \emph{from both
extremities of
the hypercube}.

To write an expression for $\Theta_k$ similar to~(\ref{expnk}), we
introduce the following notations:
\begin{itemize}
\item The ${L\choose k}$ nodes $k$ steps away from the
starting point are indexed by $\sigma$ and, as usual, their values are
written $x_\sigma$.
\item Similarly, $\tau$ indexes the ${L\choose k}$ nodes $k$ steps
away from
the end point and we note their values $1-y_\tau$.
\item $n_\sigma\in\{0,1,\ldots,k!\}$ is the number of open paths
from the starting point to
node~$\sigma$. (Contrary to the tree, there are several paths leading to
each node~$\sigma$.)
\item Similarly, $m_\tau$ is the number of open paths from node~$\tau
$ to
the end point.
\item $\mathbh{1}({\sigma\preceq\tau})$ indicates whether there is at
least one directed path (open or not) from node $\sigma$ to node $\tau$.
\end{itemize}
Then,
%
\begin{equation}\label{nkhyper}
\Theta_k=\sum_{|\sigma|=k} \sum
_{|\tau|=L-k}n_\sigma m_\tau\mathbh{1}({\sigma
\preceq\tau} ) (L-2k) (1-y_\tau-x_\sigma)^{L-2k-1}
\mathbh{1}(x_\sigma+y_\tau\le1),
\end{equation}
where
$\mathbh{1}({\sigma\preceq\tau})(L-2k)(1-y_\tau-x_\sigma
)^{L-2k-1}\mathbh{1}(x_\sigma+y_\tau<1)$
is the expected number of open
paths from $\sigma$ to $\tau$ given the values $x_\sigma$ and
$y_\tau$.

Our proof can be decomposed into three steps:
\begin{itemize}
\item First, we show that, as in the tree, the distribution of
$\Theta/L$ as $L\to\infty$ is the same as the distribution of
$\Theta_k/L$ as $L\to\infty$ and then $k\to\infty$.
\item Then, we show that the double sum in (\ref{nkhyper}) can be modified
(without changing the limit, of course) into a product of two sums. This
means that asymptotically $\Theta_k$ can be written as a contribution
from the
$k$ first levels (the sum on $\sigma$) times an independent contribution
from the $k$ last levels (the sum on $\tau$).
\item Finally, we show that each of these two contributions is
asymptotically identical in distribution to what we computed on the tree.
\end{itemize}

\subsection{First step: \texorpdfstring{$\Theta_k$}{$Theta_k$} and \texorpdfstring{$\Theta$}{$Theta$} have asymptotically the
same distribution}

We show in this section that, when the starting point scales with $L$
as $x=X/L$ for $X$
fixed,
%
\begin{equation}\label{1ststepa}
\lim_{L\to\infty} \frac{\Theta_k}L \mathop{\longrightarrow}_{k\to\infty}^{\mathrm{weakly}} \lim_{L\to\infty}\frac{\Theta}L.
\end{equation}
Following the same argument as on the tree, it is sufficient to show
that $\Theta_k /L$ has a weak limit (when $L\to\infty$ and then
$k\to\infty$) and
that
%
\begin{equation} \label{1ststep}
\lim_{k\to\infty}\limsup_{L\to\infty}\frac{1}{L^2}
\mathbb{E}^{{X}/L} \bigl[\operatorname{Var}(\Theta|\mathcal{F}_k)
\bigr]=0.
\end{equation}
First, remark that
%
\begin{equation}\label{diffof2ndmoments}
\mathbb{E}^x \bigl[\operatorname{Var}(\Theta|\mathcal{F}_k)
\bigr] =\mathbb{E}^x\bigl[\Theta^2\bigr]-
\mathbb{E}^x\bigl[\Theta_k^2\bigr],
\end{equation}
where we used $ \Theta_k=\mathbb{E}[\Theta|\mathcal{F}_k]$.

Second moments as in
(\ref{diffof2ndmoments}) can be written as sums over pairs of paths. For
a given path $\alpha$, we call $x^\alpha_i$ the value on the node at
step $i$
on
path $\alpha$ ($0\le i\le L$, with $x^\alpha_0=x$ and $x^\alpha
_L=1$) and
$\xi^\alpha_{i,j}$ the indicator function that path $\alpha$ is open from
steps $i$ to $j$:
%
\begin{equation}
\xi^\alpha_{i,j}=\mathbh{1}\bigl(x^\alpha_i
\le x^\alpha_{i+1}\le x^\alpha_{i+2}\le\cdots\le
x^\alpha_j\bigr).
\end{equation}
Clearly,
%
\begin{equation}
\Theta=\sum_\alpha\xi^\alpha_{0,L},
\qquad \Theta_k = \sum_\alpha\mathbb{E} \bigl[
\xi^\alpha_{0,L}|\mathcal{F}_k \bigr].
\end{equation}
We now have the following expression for the second moment:
%
\begin{equation}\label{Extheta2}
\mathbb{E}^x \bigl[\Theta^2 \bigr]=\sum
_{\alpha,\beta} \mathbb{E}^x \bigl[\xi^\alpha_{0,L}
\xi^\beta_{0,L} \bigr] =L!\sum_\alpha
\mathbb{E}^x \bigl[\xi^\alpha_{0,L}
\xi^0_{0,L} \bigr],
\end{equation}
where, by symmetry, we chose one particular arbitrary fixed path which
bears the index~0. Similarly,
%
\begin{equation}\label{Exthetak2}
\mathbb{E}^x\bigl[\Theta_k^2\bigr] =L!\sum
_\alpha\mathbb{E}^x \bigl[\mathbb{E} \bigl[
\xi^\alpha _{0,L}|\mathcal{F}_k \bigr] \mathbb{E}
\bigl[\xi^0_{0,L}|\mathcal{F}_k \bigr] \bigr].
\end{equation}
We write now $\xi_{0,L}^\alpha=\xi^\alpha_{0,k}\xi^\alpha_{k,L-k}
\xi^\alpha_{L-k,L}$. The first and last terms are $\mathcal{F}_k$-measurable,
hence
%
\begin{equation}
\mathbb{E}^x\bigl[\Theta_k^2\bigr] =L!\sum
_\alpha\mathbb{E}^x \bigl[
\xi^\alpha_{0,k}\xi^0_{0,k} \mathbb{E}
\bigl[\xi^\alpha_{k,L-k}|\mathcal{F}_k \bigr]
\mathbb{E} \bigl[\xi^0_{k,L-k}|\mathcal{F}_k
\bigr] \xi^\alpha _{L-k,L}\xi^0_{L-k,L}
\bigr].
\end{equation}
We make the same decomposition on $\xi^\alpha_{0,L}$ in (\ref{Extheta2}).
Writing
$\mathbb{E}^x[\cdot]=\mathbb{E}^x [\mathbb{E}[\cdot|\mathcal
{F}_k] ]$ and pushing out of the inner
expectation the $\mathcal{F}_k$-measurable terms, one gets
%
\begin{equation}
\mathbb{E}^x\bigl[\Theta^2\bigr] =L!\sum
_\alpha\mathbb{E}^x \bigl[\xi^\alpha_{0,k}
\xi^0_{0,k} \mathbb{E} \bigl[\xi^\alpha_{k,L-k}
\xi^0_{k,L-k}|\mathcal{F}_k \bigr]
\xi^\alpha_{L-k,L}\xi^0_{L-k,L} \bigr].
\end{equation}
Using (\ref{diffof2ndmoments}),
%
\begin{eqnarray}
&& \mathbb{E}^x \bigl[\operatorname{Var}(\Theta|\mathcal{F}_k)
\bigr]
\nonumber
\\
\label{bigdiff}
&& \quad= L!
\sum_\alpha\mathbb{E}^x \bigl[
\xi^\alpha_{0,k}\xi^0_{0,k} \bigl(
\mathbb{E} \bigl[\xi^\alpha_{k,L-k} \xi^0_{k,L-k}|
\mathcal{F}_k \bigr]\\
&&\hspace*{56pt}\qquad\qquad{}-\mathbb{E} \bigl[\xi^\alpha
_{k,L-k}|\mathcal{F}_k \bigr] \mathbb{E} \bigl[
\xi^0_{k,L-k}|\mathcal{F}_k \bigr] \bigr)
\xi^\alpha _{L-k,L}\xi^0_{L-k,L} \bigr].\nonumber
\end{eqnarray}
For a given path $\alpha$, the central term (in parenthesis) in the last
expression is a kind of covariance. Clearly, if the paths $\alpha$ and $0$
do not meet in the interval $\{k,\ldots, L-k\}$, the variables
$\xi^\alpha_{k,L-k}$ and $\xi^0_{k,L-k}$ are independent and the
covariance is zero. Therefore, \emph{we can restrict the sum over
$\alpha$ in~\textup{(\ref{bigdiff})} to the paths which cross at least once the
path~$0$ in the interval $\{k,\ldots, L-k\}$}. With this modified sum,
we can now find an upper bound on (\ref{diffof2ndmoments}). Dropping all
the negative terms and undoing the decomposition of $\xi_{0,L}^\alpha$ into
three parts, we get
%
\begin{equation}\label{sumprimed}
\mathbb{E}^x \bigl[\operatorname{Var}(\Theta|\mathcal{F}_k)
\bigr]\le L! {\sum_\alpha}'
\mathbb{E}^x \bigl[\xi_{0,L}^\alpha
\xi_{0,L}^0 \bigr],
\end{equation}
where the prime on the sum indicates that $\alpha$ runs only over all the
paths that meet path~0 at least once in $\{k,\ldots, L-k\}$.

We now bound (\ref{sumprimed}). Let $I_{p,q}$ be the set of all the paths
such that
\begin{itemize}
\item the $p+1$ first nodes (including the origin) are the same as for
path~0 (in other words, the first $p$ steps are the same as in path~0),
\item the next $L-p-q-1$ nodes are different from those of path~0,
\item the next $q+1$ nodes (thus including the end point) are the same as
for path~0.
\end{itemize}
By construction, for $p<k$ and $q<k$, a path in $I_{p,q}$ do not meet
path~0 in $\{k,\ldots, L-k\}$. Therefore
%
\begin{equation}\label{sumIpq}
\mathbb{E}^x \bigl[\operatorname{Var}(\Theta|\mathcal{F}_k)
\bigr] \le L!{\sum_{\alpha}}\mathbb{E}^x
\bigl[\xi^\alpha_{0,L}\xi ^0_{0,L} \bigr]
-L!\sum_{p=0}^{k-1}\sum
_{q=0}^{k-1}\sum_{\alpha\in
I_{p,q}}
\mathbb{E}^x \bigl[\xi^\alpha_{0,L}
\xi^0_{0,L} \bigr].
\end{equation}
Notice that the first sum is not primed; it runs over all the $L!$
possible paths $\alpha$. The inequality holds because in
(\ref{sumprimed}) we were summing over all the paths except \emph
{all} of
those not crossing path~0 in $\{k,\ldots, L-k\}$, while in (\ref{sumIpq}) we sum
over all the paths except \emph{some} of those not crossing path~0 in
$\{k,\ldots, L-k\}$.

Heuristically, the reason for which this bound is sufficient can be
read in Hegarty--Martinsson's paper \cite{hegarty-martinsson}: they
showed that the second moment of $\Theta$ is dominated by all the
pairs of paths that follow each other for some time, diverge close to
the start point, travel separately for most of the hypercube, meet
again close to the end point and then stick together. The second term
in the right-hand side of (\ref{sumIpq}) for large $k$ are precisely
those paths, and their contribution therefore sums up to the whole
second moment.

The first term in (\ref{sumIpq}) is simply $\mathbb{E}^x[\Theta^2]$, see
(\ref{Extheta2}). From \cite{hegarty-martinsson} and (\ref{varhyper}), it has the
following large $L$ limit
%
\begin{equation}\label{1stterm}
\lim_{L\to\infty}\frac{1}{L^2}\mathbb{E}^{{X}/L}\bigl[
\Theta^2\bigr]= \lim_{L\to\infty}\frac{1}{L^2} L!{\sum
_{\alpha}}\mathbb{E}^{{X}/L} \bigl[
\xi^\alpha_{0,L}\xi ^0_{0,L} \bigr] =4
\mathrm{e}^{-2X}.
\end{equation}
We now focus on the second term. For a path $\alpha$ in $I_{p,q}$, a
direct calculation shows that
%
\begin{equation}
\mathbb{E}^x \bigl[\xi^\alpha_{0,L}
\xi^0_{0,L} \bigr]=\frac
{(1-x)^{2L-p-q-2}}{(2L-p-q-2)!} \pmatrix{2L-2p-2q-2
\vspace*{2pt}\cr
L-p-q-1}.
\end{equation}
Indeed, excluding the starting and end points, there are $2(L-1)-p-q$
total different nodes in the paths $\alpha$ and 0. All these nodes must
be larger than $x$, hence the $(1-x)^{2L-p-q-2}$ term. This is however
not sufficient because the values on the nodes must be correctly ordered.
Out of the $(2L-p-q-2)!$ possible orderings (see the denominator), the
only good ones are those such that the $p$ smallest terms be in their
correct order
in the first shared segment (only one choice), the $q$ largest terms be well
ordered in the second shared segment (only one choice), and the
$2L-2p-2q-2$ remaining terms be separated into two well ordered blocks of
$L-p-q-1$ terms, one for each path; the only freedom is to choose which
terms go to path $\alpha$ and which to path~0, hence the binomial
coefficient.

To count the number of paths $\alpha$ in $I_{p,q}$, observe that the
$p$ first steps and the $q$ last steps of $\alpha$ are fixed and one
only has to choose the order in which the $L-p-q$ intermediary steps
are taken. We thus write $B(L-p-q)$ for the cardinal of $I_{p,q}$ where
$B(n)$ is the
number of permutations of $n$ elements such that for any $m$ in
$\{1,\ldots,n-1\}$ the image of $\{1,\ldots,m\}$ through the permutation
is not $\{1,\ldots,m\}$ (this ensures that $\alpha$ does not meet the
distinguished path at the $m$th intermediary step). The function
$B(n)$ is defined in \cite{oeis} and the first terms of the sequence are
$B(1)=1$, $B(2)=1$, $B(3)=3$, $B(4)=13$, $B(5)=71$. Hegarty--Martinsson
\cite{hegarty-martinsson} call this $T(n,1)$ and show
(Proposition 2.5) that $B(n)\sim n!$. Then
\begin{eqnarray*}
&& L!\sum_{\alpha\in
I_{p,q}}\mathbb{E}^x \bigl[
\xi^\alpha_{0,L}\xi^0_{0,L} \bigr] \\
&& \quad=
\frac{L!}{(L-p-q-1)!}\times\frac{(2L-2p-2q-2)!}{(2L-p-q-2)!} \times\frac{B(L-p-q)}{(L-p-q-1)!}
\times(1-x)^{2L-p-q-2}.
\end{eqnarray*}
Take $x=X/L$ and $L$ large with $p$ and $q$ fixed. The terms on the
right-hand side are
respectively equivalent to $L^{p+q+1}$, $(2L)^{-p-q}$, $L$ and $\mathrm{e}^{-2X}$,
so that
%
\begin{equation}
\lim_{L\to\infty}\frac{1}{L^2}L!\sum
_{\alpha\in
I_{p,q}}\mathbb{E}^{{X}/L} \bigl[\xi^\alpha_{0,L}
\xi ^0_{0,L} \bigr] = \frac{\mathrm{e}^{-2X}}{2^{p+q}},
\end{equation}
and
%
\begin{equation}\label{2ndterm}
\lim_{L\to\infty}\frac{1}{L^2}L!\sum
_{p=0}^{k-1}\sum_{q=0}^{k-1}
\sum_{\alpha\in
I_{p,q}}\mathbb{E}^{{X}/L} \bigl[
\xi^\alpha_{0,L}\xi ^0_{0,L} \bigr] =
\mathrm{e}^{-2X}4\bigl(1-2^{-k+1}+4^{-k}\bigr).
\end{equation}
Using (\ref{1stterm}) and (\ref{2ndterm}) in (\ref{sumIpq}), we finally
get
%
\begin{equation}
\limsup_{L\to\infty}\frac{1}{L^2} \mathbb{E}^{{X}/L}
\bigl[\operatorname{Var}(\Theta|\mathcal{F}_k) \bigr] \le
\frac{8 \mathrm{e}^{-2X}}{2^k},
\end{equation}
from which one gets (\ref{1ststep}). We are now going to show that
$\Theta_k/L$ has a weak limit (which we compute) thus yielding the
weak limit of $\Theta/L$ by (\ref{1ststepa}).

\subsection{Second step: Separating the start and the end of the hypercube}

We go back to the expression $\Theta_k$ given in (\ref{nkhyper}):
%
\begin{equation}
\Theta_k=\sum_{|\sigma|=k} \sum
_{|\tau|=L-k}n_\sigma m_\tau\mathbh{1}({\sigma
\preceq\tau} ) (L-2k) (1-y_\tau-x_\sigma)^{L-2k-1}
\mathbh{1}(x_\sigma+y_\tau\le1),
\end{equation}
and we introduce the following slightly different quantity
%
\begin{equation}
\tilde\Theta_k=\sum_{|\sigma|=k} \sum
_{|\tau|=L-k}n_\sigma m_\tau L(1-y_\tau-x_\sigma+x_\sigma
y_\tau)^{L-2k-1}.
\end{equation}
(Compared to $\Theta_k$, this one has no
$\mathbh{1}({\sigma\preceq\tau} )$, no $\mathbh{1}(x_\sigma
+y_\tau\le1)$, a
factor $L$ instead of $L-2k$ and an extra $x_\sigma y_\tau$ in the power.)
Clearly, $\Theta_k\le\tilde\Theta_k$. Furthermore, we know that
$\mathbb{E}^x[\Theta_k]=\mathbb{E}^x[\Theta]=L(1-x)^{L-1}$ so that
%
\begin{equation}
\lim_{L\to\infty}\mathbb{E}^{{X}/L} \biggl[\frac{\Theta
_k}L
\biggr]= \mathrm{e}^{-X}.
\end{equation}
Let us compute the same expectation for $\tilde\Theta_k$. Using
%
\begin{equation}
\mathbb{E}^x(n_\sigma|x_\sigma)=k(x_\sigma-x)^{k-1}
\mathbh {1}(x_\sigma\ge x),\qquad \mathbb{E}^x(m_\tau|y_\tau)=k(y_\tau)^{k-1},
\end{equation}
one gets
\begin{eqnarray}
\mathbb{E}^x \biggl[\frac{\tilde\Theta_k}L \biggr] &=&\pmatrix{L
\cr
k}
\pmatrix{L
\cr
k} \int_x^1\mathrm{d}x_\sigma
\int_0^1\mathrm{d}y_\tau
k(x_\sigma-x)^{k-1} k(y_\tau)^{k-1}(1-y_\tau-x_\sigma+x_\sigma
y_\tau)^{L-2k-1}
\nonumber
\\[-8pt]
\\[-8pt]
\nonumber
&=& \biggl[\frac{L!(L-2k-1)!}{(L-k)!(L-k-1)!} \biggr]^2(1-x)^{L-k-1},
\end{eqnarray}
so that
%
\begin{equation}\label{eq724}
\lim_{L\to\infty}\mathbb{E}^{{X}/L} \biggl[\frac{\tilde\Theta
_k}L
\biggr]= \mathrm{e}^{-X}.
\end{equation}
Finally, $\tilde\Theta_k/L -\Theta_k/L$ is a non-negative random
variable with an expectation going to zero; it thus converges
to zero in probability. Therefore, in the $L\to\infty$ limit by Slutsky's
theorem, $\tilde\Theta_k/L$ and $\Theta_k/L$ have the same
distribution as
soon as one of the limits exists.

It now simply remains to notice that
%
\begin{equation}\label{decompose}
\frac{\tilde\Theta_k}L= \biggl(\sum_{|\sigma|=k}n_\sigma
(1-x_\sigma)^{L-2k-1} \biggr) \biggl( \sum
_{|\tau|=L-k} m_\tau(1-y_\tau)^{L-2k-1}
\biggr),
\end{equation}
which means that $\tilde\Theta_k/L$ can be written has a contribution
coming from
the $k$ first steps of the hypercube times an independent contribution
coming from the $k$ last steps. The contribution from the start depends
on the value~$x$ of the origin. By symmetry, the contribution from the
end has the same law as the contribution from the start with~$x=0$.

\subsection{Third step: The start of the hypercube is like a tree}

We now focus on the first term in (\ref{decompose}):
%
\begin{equation}
\phi_k=\sum_{|\sigma|=k}n_\sigma(1-x_\sigma)^{L-2k-1}.
\end{equation}
First, notice that from (\ref{eq724}) and (\ref{decompose}) one has
\begin{equation}
\lim_{L\to\infty} \mathbb{E}^{X/L}(\phi_k) \mathbb{E}^{0}(\phi_k) = \mathrm{e}^{-X}
\end{equation}
because the sum over $\tau$ in (\ref{decompose}) is by symmetry equal in law to $\phi_k$ with a starting point equal to~0.
By taking $X=0$, this implies that
\begin{equation}
\lim_{L\to\infty} \mathbb{E}^{X/L}(\phi_k) = \mathrm{e}^{-X}.
\end{equation}
The goal is to show that for a starting point $x=X/L$, in the large $L$
limit then in the large $k$ limit, this $\phi_k$ converges weakly to
$\mathrm{e}^{-X}$ times an exponential distribution. Our strategy is to compare
$\phi_k$ (defined on the first $k$ levels of the hypercube) to the
$\Theta_k/L$ of the tree by showing that in the $L\to\infty$ limit the
two quantities have the same generating function.

The difficulty, of course, is that one cannot write directly a recursion
on the generating function of $\phi_k$ as we did on the tree because the
paths after the first step are not independent. To overcome this, we
introduce another quantity $\tilde\phi_k(b)$ which is (in a sense) nearly
equal to $\phi_k$:
%
\begin{equation}\label{deftildephi}
\tilde\phi_k(b)=\sum_{|\sigma|=k}\tilde
n_\sigma(b) (1-x_\sigma)^L,
\end{equation}
where we will shortly explain the meaning of the parameter~$b$ and give the
definition of $\tilde n_\sigma(b)$. For now, let us just say that
$\tilde n_\sigma(b)\le n_\sigma$; in other words, we discard some
open paths when computing $\tilde\phi_k(b)$. It is clear that
%
\begin{equation}\label{smaller}
\tilde\phi_k(b)\le\phi_k
\end{equation}
and we will choose $\tilde n_\sigma(b)$ in such a way that
%
\begin{equation}\label{limeq}
\lim_{L\to\infty}\mathbb{E}^{{X}/L} \bigl[\tilde
\phi_k(b) \bigr] =\lim_{L\to\infty}\mathbb{E}^{{X}/L}
[\phi_k ] =\mathrm{e}^{-X}.
\end{equation}
With the same argument as before, (\ref{smaller}) and (\ref{limeq}) will
be sufficient to conclude that if $\lim_{L\to\infty}\tilde\phi_k(b)$
exists (we will show it is the case), then $\lim_{L\to\infty}\phi_k$
exists as well and has the same distribution. Then, we will be able to
write a recursion for the generating function of $\tilde\phi_k(b)$
and solve
it in the $L\to\infty$ limit.

It has been pointed out to us by an anonymous referee that an
alternative way
to obtain convergence of $\phi_k$ is to use the objective method as in
Aldous--Steele~\cite{aldous-steele} and prove that the rescaled
weighted hypercube
$\{ x_\sigma L, \sigma\in\{ 0, 1\}^L\}$ converges weakly to the
so-called Poisson
Weighted Infinite Tree. This will make the Poisson cascade representation
in Section~\ref{seccascade} more intuitive.

Let us recall the following standard representation
of the hypercube: to each node of the hypercube, we associate a different
binary word with $L$ bits (digits) in such a way that the starting
point is
$(0,0,\ldots,0)$, the end point is $(1,1,\ldots,1)$ and making a step is
changing a single zero into a one. A node~$\sigma$ at level $k$ has a label
with exactly $k$ ones.

We can now define $b$ and $\tilde n_\sigma(b)$. The parameter $b$ is a set
of forbidden bits. Any path going through any bit in $b$ is automatically
discarded. In other words, $\tilde n_\sigma(b)=0$ if $\sigma$ has any bit
equal to~1 which is in $b$. The parameter $\tilde n_\sigma(b)$ is 1 or 0,
depending on whether there is an ``interesting'' path or not to $\sigma$. An
interesting path is defined recursively in the following way:
\begin{itemize}
\item From the origin, we consider which nodes amongst the $L-|b|$
reachable first level nodes have a value which is smaller than $(\ln L)/L$;
these are the ``interesting'' nodes at first level, and only the paths
going through these interesting nodes are deemed interesting and are
counted in $\tilde n_\sigma$.
\item Let $b'$ be the bits corresponding to all the interesting nodes
at first level. After the first step, these $b'$ bits are now forbidden for
all interesting paths.
\item Given the forbidden bits, the region of the hypercube reachable
from each interesting node at first level is a sub-hypercube of dimension
$L-|b|-|b'|$. All these hypercubes are non-overlapping. The
construction of the interesting paths from each first level interesting node
is now done recursively in the same way on each corresponding sub-hypercube.
\end{itemize}
Notice that by construction $\tilde n_\sigma(b)=0$ if $x_\sigma>(\ln L)/L$.
This is a small price to pay as we expect that only the $x_\sigma$ of order
$1/L$ contribute. Furthermore, at each step we exclude $\mathcal
{O}(\ln L)$ bits.
For each open paths, at step $k$, there will therefore be $k\mathcal
{O}(\ln L)$
forbidden bits. This is very small compared to $L$ and will become
negligible in the large $L$ limit.

The definition of $\tilde n_\sigma(b)$ leads directly to a recursion on
$\tilde\phi_k(b)$:
%
\begin{equation}\label{recurse}
\tilde\phi_k(b,\mbox{starting point}=x)=\sum
_{\rho\in b'}\mathbh{1}(x\le x_\rho) \tilde
\phi_{k-1}^{(\rho)}\bigl(b\cup b',\mbox{starting
point}=x_\rho\bigr),
\end{equation}
where $b'$ is the (random) set of interesting first level nodes, those
with a value smaller than $(\ln L)/L$ which avoid the $b$ forbidden
bits. \emph{Given} $b'$, for each bit $\rho\in b'$,
$\tilde\phi_{k-1}^{(\rho)}$ is an \emph{independent} copy of the variable
defined in (\ref{deftildephi}) with a different starting point. The
recursion is initialized by
%
\begin{equation}\label{startrecurse}
\tilde\phi_0(b)=(1-x)^L,
\end{equation}
which is non-random and independent of $b$.

Before computing the expectation and the generating function, remark that
the distribution of $\tilde\phi_k(b)$ depends only on the number
$|b|$ of
forbidden bits, not on the bits themselves. We will abuse this remark and
consider from now on that in the expression $\mathbb{E}^x[\tilde\phi
_k(b)]$, the
parameter $b$ is actually the \emph{number} of forbidden bits.

Let us now compute the expectation of $\tilde\phi_k(b)$. The
distribution of
the number $b'$ of interesting nodes is binomial and we call $p(b')$ its
law:
%
\begin{equation}
p\bigl(b'\bigr)=\pmatrix{L-b
\cr
b'} \biggl(
\frac{\ln L}L \biggr)^{b'} \biggl(1-\frac{\ln
L}L
\biggr)^{L-b-b'}.
\end{equation}
Then from (\ref{recurse})
%
\begin{equation}
\mathbb{E}^x \bigl[\tilde\phi_k(b) \bigr]=\sum
_{b'=0}^{L-b}p\bigl(b'\bigr) \times
b'\int_x^{({\ln L})/L}\frac{L \,\mathrm{d}y}{\ln
L}
\mathbb{E}^y \bigl[\tilde\phi_{k-1}\bigl(b+ b'
\bigr) \bigr].
\end{equation}

We will show by recurrence that the dependence in $b$ can be written as
%
\begin{equation} \label{soluceexp}
\mathbb{E}^x \bigl[\tilde\phi_k(b) \bigr]=
\frac{(L-b)!}{
(L-b-k)! L^k} \psi_k(x,L).
\end{equation}
It is obvious from (\ref{startrecurse}) that this works for $k=0$. Assume
that it works at level $k-1$. Then
%
\begin{equation}
\mathbb{E}^x \bigl[\tilde\phi_k(b) \bigr]=
\frac{1}{L^{k-1}}\sum_{b'=0}^{L-b}p
\bigl(b'\bigr)\frac{(L-b-b')!}{
(L-b-b'-k+1)!}b'\int
_x^{({\ln L})/L}\frac{L \,\mathrm{d}y}{\ln
L}\psi_{k-1}(y,L).
\end{equation}
The sum on $b'$ decouples from  the integral and can be computed; one
finds
%
\begin{equation}
\sum_{b'=0}^{L-b}p\bigl(b'
\bigr)\frac{(L-b-b')!}{
(L-b-b'-k+1)!}b' 
=
\frac{(L-b)!}{(L-b-k)!}\frac{\ln L}L \biggl(1-\frac{\ln
L}L
\biggr)^{k-1}
\end{equation}
and one recovers (\ref{soluceexp}) with
%
\begin{equation}
\psi_k(x,L)= \biggl(1-\frac{\ln L}L \biggr)^{k-1}\int
_x^{({\ln L})/L}L \,\mathrm{d} y \psi_{k-1}(y,L)
\end{equation}
or
%
\begin{equation}\label{recpsi}
\psi_k \biggl(\frac{X} L,L \biggr)= \biggl(1-
\frac{\ln L}L \biggr)^{k-1}\int_X^{{\ln L}}
\mathrm{d}Y \psi_{k-1} \biggl(\frac{Y} L,L \biggr).
\end{equation}
From here and $\psi_0(x,L)=(1-x)^L$, it is straightforward to show by
recurrence that
%
\begin{equation}
\psi_k(X/L,L)\le \mathrm{e}^{-X}.
\end{equation}
Then, with this bound and the dominated convergence theorem, the limit of
the integral in~(\ref{recpsi}) is the integral of the limit and one shows
by another straightforward recurrence that $\lim_{L\to\infty}
\psi_k(X/L,L)=\mathrm{e}^{-X}$.

Going back to (\ref{soluceexp}), one then gets for any function $b(L)$
such that
$b(L)=\mathrm{o}(L)$
%
\begin{equation}\label{boundE}
\mathbb{E}^{{X}/L} \bigl[\tilde\phi_k(b) \bigr]\le
\mathrm{e}^{-X}, \qquad \lim_{L\to\infty} \mathbb{E}^{{X}/L}
\bigl[\tilde\phi_k \bigl(b(L) \bigr) \bigr]= \mathrm{e}^{-X}.
\end{equation}
This completes the proof that $\tilde\phi_k(b)$ and $\phi_k$ have
the same
distribution in the $L\to\infty$ limit if that limit exists.

We now compute
the distribution of $\tilde\phi_k(b)$ by writing a generating
function. For
$\mu\ge0$, let
%
\begin{equation}\label{defgen}
G_k(\mu,x,L,b)=\mathbb{E}^x \bigl[\exp \bigl(-\mu\tilde
\phi _k(b) \bigr) \bigr].
\end{equation}
(Here again, we consider that the parameter $b$ of $G_k$ is a number.)
From (\ref{recurse}),
%
\begin{eqnarray}
G_k(\mu,x,L,b) &=&\sum_{b'=0}^{L-b}p
\bigl(b'\bigr) \biggl[\frac{L}{\ln
L} \biggl(x+\int_x^{({\ln L})/L}\mathrm{d}y G_{k-1}\bigl(\mu,y,
L,b+b'\bigr) \biggr) \biggr]^{b'}
\nonumber
\\[-8pt]
\\[-8pt]
\nonumber
& =&\sum_{b'=0}^{L-b}p\bigl(b'
\bigr) \biggl[1-\frac{L}{\ln L}\int_x^{({\ln L})/L}
\mathrm{d}y \bigl[1-G_{k-1}\bigl(\mu,y, L,b+b'\bigr)
\bigr] \biggr]^{b'}.\hspace*{18pt}
\end{eqnarray}
So
%
\begin{equation} \label{recurG}
G_k \biggl(\mu,\frac{X} L,L,b \biggr)= \sum
_{b'=0}^{L-b}p\bigl(b'\bigr) \biggl[1-
\frac{1}{\ln L}\int_X^{\ln L}\mathrm{d} Y
\biggl[1-G_{k-1} \biggl(\mu,\frac{Y} L, L,b+b'
\biggr) \biggr] \biggr]^{b'}.
\end{equation}
If the $G_{k-1}(\cdots)$ on the right-hand side did not depend on $b'$,
one could compute exactly the sum on $b'$ using Newton's binomial
formula. We will write bounds on
$G_{k-1}$ using quantities that do not depend on $b'$
and compute this sum.

To do this, remark that $G_k$ is an increasing function of $b$. Indeed, as
we forbid more bits ($b$ increases), we close more open paths,
$\tilde\phi_k(b)$ decreases (or remains constant) and, from (\ref{defgen}),
$G_k$ increases.

Therefore, a lower bound is easy: $G_{k-1}(\mu,Y/L,L,b+b')\ge
G_{k-1}(\mu,Y/L,L,b)$ and
%
\begin{equation}\label{lowerbound}
G_k \biggl(\mu,\frac{X} L,L,b \biggr) \ge \biggl[1-
\frac{1} L\int_X^{\ln L}\mathrm{d}Y
\biggl[1-G_{k-1} \biggl(\mu,\frac{Y} L, L,b \biggr) \biggr]
\biggr]^{L-b}.
\end{equation}

To obtain an upper bound, we use the fact that according to $p$, the
probability that $b'$ is larger than $\ln^2L$ is very small. Then, in
(\ref{recurG}), we cut the sum over $b'$ into two contributions. In
the first part $b'$
runs from 0 to $\lfloor\ln^2 L\rfloor$ and in the second part it runs
from $\lfloor\ln^2 L\rfloor+1$ to $L-b$. In the first part, we write
$G_{k-1}(\mu,Y/L,L,b+b')\le G_{k-1}(\mu,Y/L,L,b+\lfloor\ln^2
L\rfloor)$ and extend
again the sum to $L-b$. In the second part, we write that the term
multiplying $p(b')$ is smaller than~1. Hence,
%
\begin{eqnarray}
G_k \biggl(\mu,\frac{X} L,L,b \biggr) &\le & \biggl[1-
\frac{1} L\int_X^{\ln L}\mathrm{d}Y
\biggl[1-G_{k-1} \biggl(\mu,\frac{Y} L, L,b+\bigl\lfloor
\ln^2 L\bigr\rfloor \biggr) \biggr] \biggr]^{L-b}
\nonumber
\\[-8pt]
\label{upperbound}
\\[-8pt]
\nonumber
&&{}+\sum_{b'=\lfloor\ln^2L\rfloor+1}^{L-b} p\bigl(b'
\bigr).
\end{eqnarray}
The remaining sum is of course the probability that $b'$ is larger than
$\ln^2 L$, which is vanishingly small as $b'$ is binomial of average
and of variance smaller than $\ln L$.

We can now show that $G_k(\mu,X/L,L,b)$ has a large $L$ limit by recurrence.
More precisely, we will show that for any function $b(L)$ which is
a $\mathrm{o}(L)$,
%
\begin{equation}\label{Recur}
\tilde{G}_k(\mu,X):=\lim_{L\to\infty} G_k
\biggl(\mu,\frac{X} L, L, b(L) \biggr)
\end{equation}
exists and is independent of $b(L)$.

This is obvious for $k=0$ as $G_0(\mu,x,L,b)=\exp[-\mu(1-x)^L]$, so that
%
\begin{equation}\label{inittildeG}
\tilde{G}_0(\mu,X)=\exp \bigl[-\mu \mathrm{e}^{-X}
\bigr].
\end{equation}
Suppose that (\ref{Recur}) holds up to level~$k-1$. Then for any function
$b(L)=\mathrm{o}(L)$, the function $b(L)+\lfloor\ln^2 L\rfloor$ is also an $\mathrm{o}(L)$.
We know from (\ref{defgen}) and (\ref{boundE}) that $G_k(\mu
,X/L,L,b)\ge
1 -\mu\mathbb{E}^{{X}/L}[\tilde\phi_k(b)]\ge1-\mu \mathrm{e}^{-X}$, so
that we
can use the dominated convergence theorem and obtain
%
\begin{equation}
\lim_{L\to\infty} \int_X^{\ln L}
\mathrm{d}Y \biggl[1-G_{k-1} \biggl(\mu,\frac{Y} L, L,\mathrm{o}(L)
\biggr) \biggr]=\int_X^\infty\mathrm{d}Y \bigl[1-
\tilde G_{k-1}(\mu,Y) \bigr].
\end{equation}
It is then straightforward from (\ref{lowerbound}) and (\ref{upperbound})
to see that (\ref{Recur}) holds at level $k$ and that
%
\begin{equation}\label{rectildeG}
\tilde{G}_k(\mu,X)=\exp \biggl[-\int_X^\infty
\mathrm{d}Y \bigl[1-\tilde G_{k-1}(\mu,Y) \bigr] \biggr].
\end{equation}
Equations (\ref{inittildeG}) and (\ref{rectildeG}) are the same as
(\ref{limG2}), which completes the proof.

\section*{Acknowledgements}

We thank an anonymous referee who pointed out the relevance of \cite
{aldous-steele} to our work.

%

\printhistory
\end{document}